\documentclass{article}
\usepackage[english]{babel}
\usepackage[normalem]{ulem}
\usepackage{graphicx}
\usepackage{dcolumn}
\usepackage{bm}
\usepackage{cancel}
\usepackage{amssymb, amsmath}
\usepackage{multirow}
\usepackage{mathbbol} 
\usepackage{latexsym}
\usepackage{xspace,bm}
\usepackage{hyperref}
\usepackage{cleveref}
\usepackage{dsfont}
\usepackage{xcolor}
\usepackage{authblk}
\usepackage[letterpaper,top=2cm,bottom=2cm,left=2cm,right=2cm,marginparwidth=1.75cm]{geometry}
\usepackage{enumitem}
\setlist{topsep=1pt, partopsep=0pt, itemsep=0pt, parsep=1pt}
\usepackage{todonotes}

\newcommand{\E}{\mathbb{E}}
\newcommand{\Lcoarse}{L_\mathrm{c}}
\newcommand{\MV}{\mathcal{V}}
\newcommand{\MW}{\mathcal{W}}
\newcommand{\td}{\text{d}}

\title{Trajectory learning for ensemble forecasts via the continuous ranked probability score: a Lorenz '96 case study}
\author[1]{Sagy R. Ephrati\footnote{(sagy@chalmers.se)}}
\affil[1]{Department of Mathematical Sciences, Chalmers University of Technology and University of Gothenburg, 412 96 Gothenburg, Sweden}
\author[2]{James Woodfield\footnote{(j.woodfield18@imperial.ac.uk)}}
\affil[2]{Department of Mathematics, Imperial College London, South Kensington Campus, London, SW7 2AZ, United Kingdom}

\date{August 2025}

\begin{document}

\maketitle
\begin{abstract}
    This paper demonstrates the feasibility of trajectory learning for ensemble forecasts by employing the continuous ranked probability score (CRPS) as a loss function. Using the two-scale Lorenz '96 system as a case study, we develop and train both additive and multiplicative stochastic parametrizations to generate ensemble predictions. Results indicate that CRPS-based trajectory learning produces parametrizations that are both accurate and sharp, while the self-spread term in the loss function regularises training. The resulting parametrizations are straightforward to calibrate and outperform derivative-fitting-based parametrizations in short-term forecasts. This approach is particularly promising for data assimilation applications due to its accuracy over short lead times.
\end{abstract}

\section{Introduction}

Addressing uncertainty stemming from unresolved physical effects or uncertain measurements constitutes a core challenge in geophysical fluid dynamics. Incomplete representations of physical processes, finite truncation of computational models used for dynamic simulations, and uncertain initial conditions all contribute to numerical prediction uncertainty. Stochastic parametrizations provide a means to simultaneously compensate for model error and predict the inherent uncertainty via ensemble forecasts, based on available observational data. The modeling objective thus shifts to calibrating a noise model or parametrization. In this work, we explore the calibration of stochastic models for ensemble forecasting via \textit{trajectory learning} based on the continuous ranked probability score (CRPS).

Stochastic methods have been effectively utilised in numerical weather prediction for several decades \cite{hasselmann1976stochastic, buizza1999stochastic}. In recent years, data-driven stochastic models have been actively developed. The two-scale Lorenz '96 system \cite{lorenz1996predictability}, a simplified representation of atmospheric dynamics with interactions between slow and fast scales, commonly serves as an initial test for novel calibration strategies \cite{balwada2024learning}. Neglecting the fast variables is analogous to coarsening in fluid simulations, and hence introduces error and uncertainty in the evolution of the slow variables. Consequently, data from fast scales can inform stochastic parametrizations. For instance, \cite{wilks2005effects} and \cite{arnold2013stochastic} employed explicitly computed sub-grid scale data to tune parametrizations, fitting them as polynomials with a stochastic component for residuals, modeled as correlated or uncorrelated noise. Ensemble predictions with stochastic forcing were found to outperform perturbed parameter ensembles, with correlated additive noise yielding reliable forecasts. \cite{crommelin2008subgrid} modeled the sub-grid effects using Markov chains conditional on the resolved variables, achieving good agreement in the autocorrelation of the resolved variables and in their long-time distribution. Generative adversarial networks were used in \cite{gagne2020machine} to construct maps between resolved variables and sub-grid scale forcings, producing accurate spatiotemporal correlations.

The development of data-driven stochastic models extends beyond low-dimensional dynamical systems, with significant theoretical and methodological advancements made in fluid models. Examples include approaches based on statistical turbulence theory \cite{frederiksen2024statistical} or the assumption of spatiotemporal coherent error patterns \cite{buizza1999stochastic}.
Modern approaches include developing structure-preserving stochastic models that exploit the geometry of geophysical fluid-dynamical models to construct specific multiplicative noise. For example, the framework of Stochastic Advection by Lie Transport (SALT) \cite{holm2015variational} is well suited for advection-dominated flows. This approach has yielded efficient uncertainty quantification for the two-dimensional Euler equations \cite{cotter2019numerically, ephrati2023data, woodfield2023stochastic} and (thermal) quasi-geostrophic equations \cite{cotter2018modelling, holm2024comparing} in a reduced-order setting, with applications in data assimilation \cite{cotter2020data}. 

A pragmatic low-rank approach to account for uncertainty and model error in fluid-dynamical systems involves separating spatiotemporal forcing into fixed spatial profiles and corresponding time series. These time series are subsequently modeled as stochastic processes. Proper orthogonal decomposition modes (POD) modes, also referred to as empirical orthogonal functions (EOFs) in atmospheric science \cite{tennekes1972first, wilks2011statistical,hannachi2007empirical}, efficiently capture coherent structures observed in data. For example, \cite{agarwal2021comparison} compared data-driven methods for ocean models and found that stochastic processes with memory effects performed well in both short- and long-time forecasts. Global basis methods have also been employed to capture flow statistics in coarsened turbulence simulations, resulting in stochastic parametrizations that reproduce energy spectra, total energy, or enstrophy \cite{hoekstra2024reduced, ephrati2025probabilistic, ephrati2025continuous, ephrati2023rb}. 

In this work, we aim to produce a forecast with an CRPS and thereby tackle the model calibration challenge. While this challenge does not coincide with parameter estimation in stochastic differential equations (SDEs), they are related.  Methods based upon maximum likelyhood estimation \cite{pedersen1995new,pedersen1995consistency}, generalised method of moments \cite{hansen1982large}, and the efficient method of moments \cite{gallant1996moments} are among the approaches outlined in the review \cite{nielsen2000parameter} that bear similarity to the CRPS learning approach adopted in this paper. However, parameter estimation in SDEs differs from producing a skillful forecast. The CRPS approach used in this study, can be viewed as both a data assimilation technique for ensemble forecasting and as a model calibration strategy \cite{hourdin2017art}.

Compensating for model error and unresolved dynamics remains a longstanding challenge in computational fluid dynamics, with the resulting models generally referred to as `closure models'. Machine learning has proven to be an efficient approach for model calibration (see \cite{sanderse2024scientific} and references therein). A robust method is \textit{trajectory fitting}, where the discrepancy between predicted solution trajectories and reference trajectories is minimised \cite{melchers2023comparison}. However, many learned closure models are deterministic. Learning stochastic models is a logical next step necessary for probing prediction uncertainty via ensemble methods \cite{christensen2024machine}, and currently forms an active research area. A particular inspiration for the present study is the recent work by \cite{lang2024aifs}, who trained a stochastic model based on the \textit{almost fair} CRPS to avoid problems with score degeneracy that can arise when all but one ensemble member has the same value as the reference. This approach is employed in the Artificial Intelligence Forecasting System (AIFS) developed at the ECMWF in test cases at two resolutions, where the trained stochastic model demonstrated improved variability compared to train deterministic models. The stochastic model was optimised for short-range forecasts but also showed good performance for longer range forecasts. Similarly, in \cite{kochkov2024neural} neural global circulation models for the atmosphere are trained using a stochastic loss based on the CRPS, giving strong practical evidence suggesting that minimising the CPRS to predict short-term weather is an effective way to tune
parameterizations for climate.

\paragraph{Contributions and paper outline}
This work investigates the feasibility of calibrating stochastic parametrizations for ensemble forecasts through trajectory learning, employing the CRPS as a loss function. This choice of loss function directly incorporates forecast reliability and skill into the training procedure. A simple linear stochastic model form is used to ensure interpretability and serves as a proof of concept for CRPS-based trajectory learning. Specifically, a multidimensional coupled Ornstein--Uhlenbeck process is adopted in additive and multiplicative settings to account for unresolved dynamics and estimate uncertainty. The presented approach and results inform future studies on stochastic parametrization in fluid models based on trajectory learning of ensemble forecasts.

The paper is structured as follows. The noise calibration methodology is introduced in Section \ref{sec:calibration_methodology}. The governing equations are provided in Section \ref{sec:governing_equations}, along with the adopted parametrizations and calibration details. The model performance is assessed in Section \ref{sec:model_performance}, after which the paper is concluded in Section \ref{sec:conclusion}.

\section{Noise calibration methodology}\label{sec:calibration_methodology}
This section introduces the noise calibration methodology. We first summarize the derivation of a general deterministic sub-grid scale measurement procedure in Section \ref{subsec:deterministic_measurements}. This is then extended to a probabilistic procedure in Section \ref{subsec:probabilistic_measurements}, where the use of the CRPS as a loss function is motivated.

\subsection{Deterministic sub-grid scale measurement procedure} \label{subsec:deterministic_measurements}
We consider two dynamical systems.
The first is the `true' dynamics evolving according to $\dot q=L(q)$.
The second is the `coarse' dynamics described by $\dot{z} = \Lcoarse(z)$.
The true dynamics is integrated in time using the `high-fidelity solver', whereas the coarse dynamics is integrated using the `low-fidelity solver'. 
Corresponding solution snapshots are referred to as high-fidelity and low-fidelity snapshots, respectively.
We introduce the solution spaces $\MV$ for the true solution and $\MW$ for the coarse solution. This formulation allows $\MV$ and $\MW$ to be either function spaces or spaces of solutions on discrete meshes. In the following, we treat these as the latter, where elements of $\MV$ are high-resolution numerical solutions and elements of $\MW$ are low-resolution numerical solutions.
We further introduce a coarsening operator $C:\MV\to\MW$, which is generally assumed to discard information and is thus not invertible.

The true dynamics will generally be either unknown or too computationally expensive for practical forecasting.
On the other hand, the coarse dynamics provides a computationally feasible alternative but suffers from model incompleteness and other sources of error.
Therefore, a forcing term will be included in the coarse dynamics to compensate for these errors. 
Generally, we refer to the discrepancy between $L_c$ and $L$ as `unresolved scales' leading to `unresolved dynamics'.

A perfect model is one where the modeled coarse dynamics equal the coarsened true dynamics. This requires defining a model correction term, or parametrization, $M(C(q), \theta)$, such that the complete coarse dynamics $L_c + M$ satisfies
\begin{equation}\label{eq:perfect_model}
    L_c(C(q))+M(C(q), \theta) = C(L(q)), 
\end{equation}
where $\theta$ is a tunable set of parameters.

It is generally not possible to find closed-form expressions for this perfect model, as it depends on properties of the governing dynamics as well as on adopted numerical methods.
Nonetheless, approximations to this model can be found using observational data. A direct approach based on \eqref{eq:perfect_model} is \textit{derivative fitting}, where the parameters $\theta$ are chosen to minimize the loss function \begin{equation}
    \mathrm{Loss}(\theta, q) = \left| L_c(C(q) + M(C(q), \theta) - C(L(q)) \right|^p.
\end{equation}
That is, the objective is to match the instantaneous coarse-grained dynamics over a given training data set. Here, $p$ determines how strongly outliers are penalised. Naturally, regularization terms can be included in the loss function without loss of generality, but they will not be further considered here.

An alternative approach is \textit{trajectory learning}, where the model aims to minimize the distance between the low-fidelity dynamics and observed coarse-grained high-fidelity dynamics over specified integration intervals. 
We denote a reference trajectory by $q(t_j),~j=1,\ldots,N_t$, generated by the dynamics $\dot q = L(q)$, and a model trajectory by $z_\theta(t_j),~j=1,\ldots,N_t$, generated by the low-fidelity dynamics including the model as $\dot z = L_c(z) + M(z, \theta)$. Then, the model parameters $\theta$ are chosen such that the loss function \begin{equation}
    \mathrm{Loss}(\theta, q) = \sum_{j=1}^{N_t} \left| z_\theta(t_j) - C(q(t_j))\right|^p
\end{equation} is minimised over trajectories within a training data set. Here, $N_t$ denotes the length of each trajectory. Trajectory fitting can thus account for long-term dynamics (temporal non-locality) and error accumulation when the trajectory length is appropriately chosen. 

A concise comparison between derivative fitting and trajectory learning is provided in \cite[Table 1]{melchers2023comparison}. In brief, the difference between derivative fitting and trajectory learning decreases when the trajectory length is shortens and the time instances $t_j$ become sufficiently close. However, when data are sparse or noisy, approximations of the instantaneous dynamics based on available measurements become less accurate and derivative fitting loses reliability. Simultaneously, derivative fitting requires time derivatives of training data, which may not always be readily available. Conversely, trajectory fitting requires either a discretised adjoint equation or a differentiable numerical solver. Furthermore, if the trajectory length is too long or the time instances are spaced too far apart, chaotic behavior and perturbation sensitivity of the dynamical system may decrease the effectiveness of trajectory learning.

\subsection{Probabilistic sub-grid scale measurement procedure} \label{subsec:probabilistic_measurements}
The assumption that the coarsening operator $C$ discards information has considerable consequences for model derivation. Namely, the preimage of a coarse configuration $z$ through $C$ is a distribution of high-fidelity solutions $\{q_i\}_i\in\MV$. Since each $q_i$ satisfies $C(q_i) = z$, it is impossible to determine which high-fidelity solution produced the particular coarse configuration $z$. Hence, there is uncertainty in the `perfect' evolution of a coarse-grained model. This motivates a probabilistic approach to sub-grid model derivation, aiming to match \textit{distributions of solutions} rather than a single solution. Stochastic modeling is well-suited for this purpose.

This intuition can be made more precise as follows. We let $\pi_\mathrm{full}$ denote the distribution of fully resolved configurations, where $\pi_\mathrm{full}^{t_0}(q)$ describes the probability of a configuration $q$ being the realization of the dynamics at time $t_0$. The true dynamics are defined via a high-fidelity flow map $\psi_{[t_0, t_1]}(q):= q + L_{[t_0, t_1]}(q)$, evolving the dynamics from $t=t_0$ to $t=t_1$ using the high-fidelity solver. 

Given a coarse configuration $z_0$ at $t_0$, we denote the set of corresponding full configurations by $\mathcal{A}=\{q\in\mathcal{V} \mid C(q)=z_0 \}$ and its distribution by $\pi_\mathcal{A}^{t_0}$. The high-fidelity flow map is deterministic, hence the evolution of this distribution under the high-fidelity dynamics can be expressed via a transition kernel $\tau(q'|q) = \delta(q' - \psi_{[t_0, t_1]}(q))$. At time $t_1$, the distribution is expressed as \begin{equation}
    \pi_\mathcal{A}^{t_1}(q') = \int_\mathcal{A}\! \delta(q' - \psi_{[t_0, t_1]}(q))\pi_\mathcal{A}^{t_0} \,\text{d}q.
\end{equation}
The configurations that make up $\pi_\mathcal{A}^{t_1}$ are elements of $\mathcal{V}$, i.e., these are fully resolved configurations. The corresponding distribution of coarse-grained fields is obtained by marginalizing, \begin{equation}
    \pi_\mathrm{CG}^{t_1}(z) = \int_\mathcal{B}\! \pi_\mathcal{A}^{t_1}(q)\,  \text{d}q,
\end{equation}
where $\mathcal{B} = \{q\in\mathcal{V} \mid C(q)=z\}$. Intuitively, $\pi_\mathrm{CG}^{t_1}(z)$ is the distribution of attainable coarse-grained fields at time $t_1$, given the coarse configuration $z_0$ at $t_0$ and the uncertainty induced by the coarsening operator $C$.

Hence, from this perspective, the sub-grid scale model and the low-fidelity solver should define a transition kernel $\tau_\mathrm{model}$ that satisfies $\tau_\mathrm{model}(z'|z) = \pi_\mathrm{CG}^{t_1}(z')$. In this work, the transition kernel is defined implicitly via the stochastic forward flow map \begin{equation}\label{eq:flow_map_stochastic}
    \varphi_{[t_0, t_1]}(z, \theta) := z + L_{c, [t_0, t_1]}(z)  + M_{[t_0, t_1]}(z, \theta),
\end{equation} 
where $L_{c, [t_0, t_1]}$ denotes low-fidelity integration from time $t_0$ to $t_1$. The sub-grid scale model is denoted by $M$ and depends on parameters $\theta$. Subsequently, an ensemble prediction obtained with the flow map \eqref{eq:flow_map_stochastic} will serve as a Monte Carlo approximation to the distribution $\pi_\mathrm{CG}^{t_1}$.

The distribution $\pi_\mathrm{CG}^{t_1}$ of coarse-grained field will generally not be known. Monte Carlo approximations can be computed using the high-fidelity flow map $\psi_{[t_0, t_1]}$, yet this map can be unknown or may be prohibitively expensive to compute. Instead, a single observed trajectory can be treated as a time-varying delta distribution to facilitate probabilistic trajectory learning. To this end, we consider an ensemble forecast $\{z_\theta^i(t_j)\}_{i=1,\ldots, M},~ j=1,\ldots, N_t$, generated by the flow map \eqref{eq:flow_map_stochastic}. Here, the superscript $i$ denotes the $i^\mathrm{th}$ ensemble member, $M$ is the ensemble size, and $N_t$ is the trajectory length. As before, we let $q(t_j), ~j=1,\ldots, N_t,$ be the reference trajectory generated by the high-fidelity dynamics. The continuous ranked probability score (CRPS) then serves as a probabilistic extension to the mean absolute error, enabling trajectory learning for ensemble predictions. The objective becomes to choose model parameters $\theta$ such that the loss function \begin{equation} \label{eq:CRPS_loss}
    \begin{split}
        \mathrm{Loss}(\theta, q) &= \sum_{j=1}^{N_t} \mathrm{CRPS}(z_\theta, C(q), t_j), \\
        \mathrm{CRPS}(z_\theta, q, t_j) &= \frac{1}{M}\sum_{i=1}^M \left|z_\theta^i(t_j) - C(q(t_j))\right| - \frac{1}{2M^2}\sum_{j=1}^M \sum_{k=1}^M\left|z_\theta^i(t_j) - z_\theta^k(t_j)\right|,
    \end{split}
\end{equation}
is minimised over a set of training trajectories. 

\begin{figure}[h]
    \centering
    \includegraphics[width=0.49\linewidth]{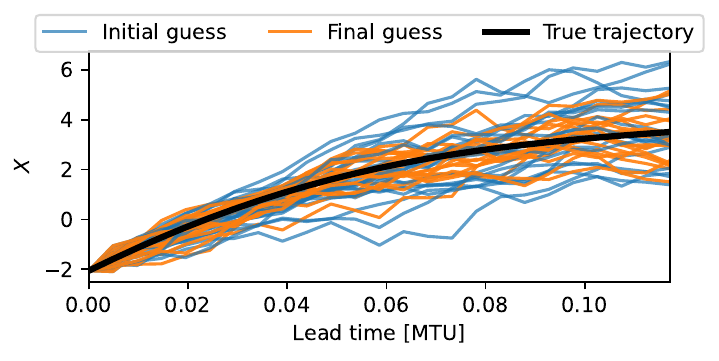}
    \includegraphics[width=0.49\linewidth]{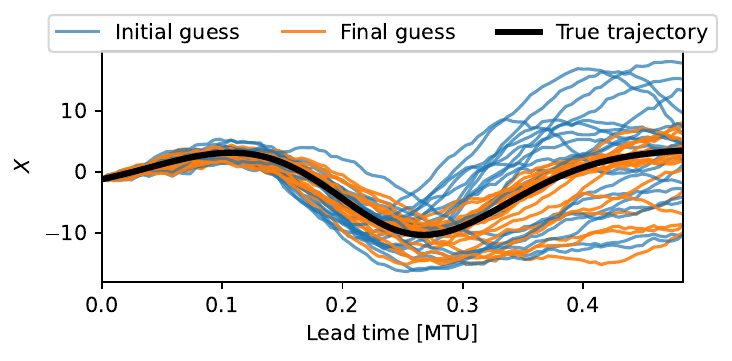}
    \caption{Example of CRPS-based trajectory learning for ensemble predictions. The true trajectory of a single variable in the L96 system serves as a reference. The trajectory length equals 25 coarse-grid time steps (left) and 100 coarse-grid time steps (right), with each ensemble consisting of 20 members. After training, the ensemble exhibits reduced mean absolute error and spread.}
    \label{fig:example_training}
\end{figure}

The first term in the CRPS represents the ensemble mean absolute error, while the second term quantifies the ensemble spread and is subtracted from the mean error. On its own, a larger spread would reduce the CRPS. However, in practice, an increase in spread usually leads to a larger mean error, which causes the CRPS to increase.
Hence, a small CRPS implies both an accurate ensemble mean absolute error and a sharp ensemble forecast (small spread). Furthermore, the CRPS is suitable for ensemble forecasts of nonlinear systems, as no assumptions are made on the distribution of the samples. As a proper scoring rule, the CRPS is minimised only when the forecast distribution matches the observed distribution. This ensures that training with the aforementioned loss function favors both ensemble accuracy and sharpness. Moreover, the absence of squared errors in the loss function makes the training less sensitive to large outliers, thereby increasing training robustness in chaotic systems. For additional details regarding the CRPS, we refer to \cite{arnold2024decompositions,gneiting2007strictly,gneiting2011comparing,hersbach2000decomposition,matheson1976scoring}.
An illustration of CRPS-based trajectory learning of ensemble forecasts is provided in Fig. \ref{fig:example_training}.

\section{Governing equations and models} \label{sec:governing_equations}
In this section, we detail the governing equations and stochastic parametrizations. We differentiate between `true' reference dynamics, given by the deterministic two-scale Lorenz '96 (L96) model, integrated at high temporal resolution. A computationally cheaper `coarsened' forward model comprises a single-scale L96 model, integrated at low temporal resolution. Thus, the latter exhibits deficiencies due to reduced temporal resolution and the omission of fast-scale dynamics. These deficiencies are parametrised by means of stochastic models, using synthetic data generated from the high-fidelity reference dynamics.

\subsection{Governing equations}

\paragraph{True dynamics}
The two-scale L96 system describes the interaction of variables across multiple spatial and temporal scales.
The large (slow) scales are denoted by $X$, and these are coupled to the small (fast) scales,  denoted by $Y$.
The governing equations read 
\begin{align}
    \frac{\td X_k}{\td t} &= -X_{k-1}(X_{k-2} -X_{k+1})-X_k + F - \frac{hc}{b}\sum_{j=J(k-1)+1}^{kJ}Y_j; &k=1,\ldots,K; \label{eq:truth1}\\
    \frac{\td Y_j}{\td t} &= -cb Y_{j+1}(Y_{j+2}-Y_{j-1})-cY_j + \frac{hc}{b}X_{\mathrm{int}[(j-1)/J]+1}; &j= 1,\ldots,JK.\label{eq:truth2}
\end{align}
Periodic boundary conditions are imposed on the variables, which translates to $X_{k+K} = X_k$ and $Y_{j+JK}=Y_j$. The two-scale L96 system (\cref{eq:truth1,eq:truth2}) is a low-dimensional deterministic representation of atmospheric dynamics at a latitude circle. The slow scale evolution of $X$ is evolved with the nonlinear term $X_{k-1}(X_{k-2} -X_{k+1})$, which mimics advection in a fluid flow. The term $-X_k$ represents a linear damping, and $F$ denotes a constant forcing, which for sufficiently large values results in chaotic dynamics. In the two-scale extension, each slow variable $X_k$ is associated with a set of fast variables $Y_{j,k}$. The fast variables in turn interact nonlinearly with themselves. This system was studied by Lorenz and Emanuel \cite{lorenz1998optimal} to determine error growth in ensemble forecasting. For a discussion of how to interpret the parameters we refer to \cite{lorenz2005designing}.

The two parameter sets considered in this study are identical to those presented by \cite{arnold2013stochastic} and are summarised in Table \ref{tab:parameters_L96}. Time is expressed in terms of model time units (MTU), where one MTU corresponds to approximately 5 atmospheric days.

\begin{table}[h]
    \centering
    \caption{Parameter values for the L96 experiments}
    \begin{tabular}{l|l|l}
       Description & Parameter & Value \\
       \hline
       Number of $X$ variables  & $K$ & 8\\
       Number of $Y$ variables per $X$ variable  & $J$ & 32\\
       Coupling constant & $h$ & 1\\
       Forcing term & $F$ & 20\\
       Spatial scale ratio & $b$ & 10\\
       Temporal scale ratio & $c$ & 4 or 10\\
    \end{tabular}
    \label{tab:parameters_L96}
\end{table}

\paragraph{Coarse dynamics}
The distinction between slow and fast variables renders the L96 system a well-suited testbed for stochastic parametrizations.
The coarse dynamics consist of only the slow variables with an additional term $M$ to parametrize the effect of the unresolved fast variables,
\begin{equation}
    \frac{\td X_k}{\td t} = -X_{k-1}(X_{k-2} -X_{k+1})-X_k + F + M(X_k, t); \quad k=1,\ldots,K. \label{eq:parametrised model}
\end{equation}
Here, $M$ denotes the parametrization of both the missing unresolved scales $Y$ and adverse numerical effects arising from integrating the coarsened system using a large time step. The parametrization will introduce stochasticity to account for uncertainty.

\paragraph{Data gathering}
A high-fidelity data set (referred to as the `truth') is generated by integrating the true dynamics using the three-stage third-order strong stability preserving (SPPRK3) time-stepping scheme \cite{shu1988efficient} with fixed time-step size $\Delta t=0.001$ MTU. The coarse dynamics are simulated with a step size $\Delta t = 0.005$ MTU. The stochastic integrals arising from the parametrizations are interpreted in the Stratonovich sense and are integrated using a stochastic SSPRK3 scheme \cite{woodfield2024strong}, which has strong order of convergence 1/2 in the stochastic setting and reverts to deterministic order 3 in the absence of noise \cite{ruemelin1982numerical}.

The high-fidelity solution is initialised by assigning to all variables an independent sample from a standard normal distribution, after which the solution is integrated for 1500 MTU to remove any transitory behavior. 
The obtained solution at this time serves as the starting point for the data measurement procedure. Subsequently, the truth is simulated for an additional 500 MTU to generate training and verification data.

\subsection{Derivative-fitting-based models}\label{subsec:derivative_fitting_models}


Sub-grid scale models can be derived from direct measurements of the sub-grid tendency \cite{arnold2013stochastic}. This procedure measures the discrepancy between the true and coarse dynamics. Consequently, models derived from this data are based on derivative fitting.

The sub-grid tendency $U$ is approximated as follows,
\begin{equation}
    U(X, t) = -X_{k-1}(X_{k-2} -X_{k+1})-X_k + F - \left(\frac{X_k(t+\Delta t)-X_k(t)}{\Delta t} \right),
\end{equation}
and is computed exclusively from the snapshots of the true dynamics. The measurements represent the difference between the evolution of the coarse dynamics and the true dynamics, with the latter approximated by an Euler forward step of size $\Delta t$. Here, $\Delta t=0.005$ in accordance with the coarse time step. 

Treating the $X$ variables as `spatial', the measurements are stored a space-time array $\mathbb{U}$. This serves as input for two distinct approaches to stochastic parametrization, referred to as the local model and the global model. Both forecast models are of the form \begin{equation}
M(X_k, t) = M_\mathrm{det}(X_k) + r(t), \end{equation} where $M_\mathrm{det}$ is a deterministic parametrization and $r(X, t)$ denotes a stochastic term, which may be space-dependent.

\subsubsection{Local model}
The first forecast model employs a third-degree least-squares polynomial fit of $U$ in terms of $X$ and a stochastic model for the residuals. That is, $M_\mathrm{det}(X_k) = b_0 + b_1X_k + b_2X_k^2 + b_3X_k^3$, as presented by \cite{arnold2013stochastic}. This parametrization is \textit{local}, as it is applied to each $X_k$ separately, and the model term does not depend on the neighboring values.

For the first and second test cases, we obtain respectively\begin{align*}
    M_\mathrm{det}(X_k) &= -2.20  +  0.575 X_k - 0.00499 X_k^2 - 0.000216 X_k^3,  \\
    M_\mathrm{det}(X_k) &= 0.324 + 1.30 X_k - 0.0128 X_k^2 - 0.00234 X_k^3. 
\end{align*}
The stochastic components of the parametrization follow from the residuals of the polynomial fit. Specifically, two additive stochastic processes are considered, based on the sample variance $S^2$ of the residuals and their autocorrelation $\phi$. The first process consists of uncorrelated, normally distributed noise, where the increments over a single time step are drawn from $\mathcal{N}(0, S^2\Delta t)$. The second process is an Ornstein--Uhlenbeck (OU) process that takes into account measured autocorrelation of the residuals. It is discretised as a zero-mean first-order autoregressive (AR(1)) process, which evolves as \begin{equation}
    r(t+\Delta t) = \phi\hspace{0.5mm} r(t) + \sqrt{\Delta tS^2(1-\phi^2)}~\varepsilon(t). \label{eq:ar1}
\end{equation}
Here, $\varepsilon(t)\sim\mathcal{N}(0, 1)$ and the scaling term ensures that the resulting samples have the same variance as the uncorrelated model. The value of $\phi$ is set to the autocorrelation of the residuals after $\Delta t=0.005$ MTU. All parameter values are summarised in Table \ref{tab:polyfit_params}.

\begin{table}[h]
    \centering
    \caption{Parameters for the stochastic processes used in the local model.}
    \begin{tabular}{l|l}
        $c=4$ & $c=10$ \\
        \hline
        $S^2 = 4.61$ & $S^2 = 4.04$ \\
        $\phi = 0.974$ & $\phi = 0.977$ 
    \end{tabular}
    \label{tab:polyfit_params}
\end{table}

\subsubsection{Global model}
The second forecast model is based on \textit{global} forcing, where the stochastic parametrization is defined as a sum of spatial profiles multiplied by scalar stochastic processes. This approach is motivated by turbulence, which is ultimately a nonlocal phenomenon. The nonlinear advection term in the L96 system ensures that the evolution of each slow variable is intrinsically linked to the state of its neighbors, thereby motivating parametrizations that utilize global basis functions.
Here, the spatial profiles are selected as POD modes to capture coherent patterns within the measurements. Alternative approaches can be used analogously, for instance by employing a Fourier basis, but these are not considered in this study.

\begin{figure}[h]
    \centering
    \includegraphics[width=0.49\linewidth]{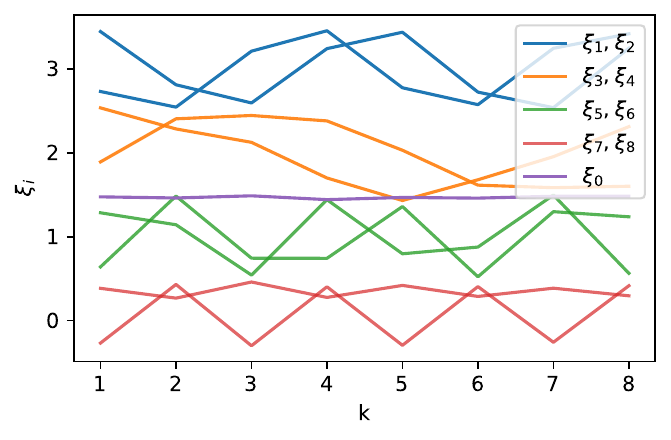}
    \includegraphics[width=0.49\linewidth]{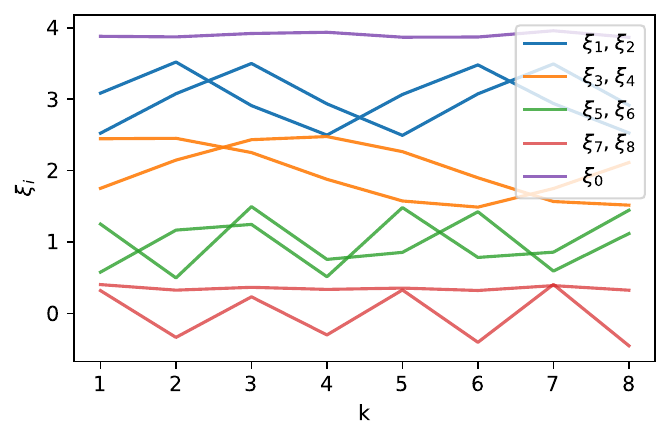}
    \caption{POD modes for $c=4$ (left) and $c=10$ (right) that serve as the basis for the global stochastic parametrization. The pairs of POD modes are shown with a vertical offset for visibility, all modes have zero mean. The mean profile $\xi_0$ does not have zero mean and is not offset.}
    \label{fig:spatial_profiles}
\end{figure}

The spatial profiles are computed from the sub-grid tendency matrix $\mathbb{U}$ in a direct manner. The time-mean is subtracted from $\mathbb{U}$, yielding the average spatial profile $\xi_0$ and an anomaly matrix $\mathbb{A}$ with zero-mean rows. The singular value decomposition $\mathbb{A}=B\Sigma V^T$ then yields spatial profiles as columns of $B$, the respective squared eigenvalues on the diagonal of $\Sigma$, and the corresponding time series as columns of $V$. The POD modes will be denoted by $\xi_i(x), i=1,\ldots,8$, and are of unit norm. The corresponding time series are denoted by $a_i$ and have unit variance. Scaling is accounted for via the eigenvalues $\lambda_i$. The decomposition of the original signal $U(X, t)$ then takes the form \begin{equation}
    U(t) = \xi_0 + \sum_{i=1}^8 \xi_i \lambda_i a_i(t).
\end{equation}
The spatial profiles are depicted in Fig. \ref{fig:spatial_profiles}, with $k$ treated as the spatial variable. It is observed that spatially coherent patterns can be extracted from the sub-grid scale data.

The time series $a_i(t)$ are replaced by independent scalar stochastic processes $r_i(t)$ to define the forecast model, analogous to the local model. This results in a sub-grid scale model of the form \begin{equation}
    M(t) = \xi_0+ \sum_{i=1}\xi_i \lambda_i r_i(t),
\end{equation}
where uncorrelated Gaussian noise, $r_i(t)\sim\mathcal{N}(0,\Delta t)$, is the first approach taken. The second approach uses correlated noise, implemented by replacing the process parameter $\phi$ in \eqref{eq:ar1} by the measured autocorrelation of $a_i(t)$. An overview of the model parameters is given in Appendix \ref{app:details}.

\subsection{CRPS-based forecast model}

The CRPS-based forecast models are defined in terms of the POD modes with coefficients determined by a learned stochastic process. Here, we consider both additive and multiplicative model forms. Generally, these models can be expressed as
\begin{equation}
    M(t) = \sum_{i}\xi_i c_i(X)\mathbf{r}_i(t),
\end{equation}
where $\mathbf{r}_i$ is the $i^\mathrm{th}$ component of coupled OU processes. The function $c_i(X)$ introduces a possibly nonlinear scaling based on the current configuration $X$.

The adopted additive noise model follows from
\begin{equation}
    \mathbf{r}(t+\Delta t) = \mathbf{r}(t) + \mathbf{A}(\boldsymbol\mu-\mathbf{r}(t)) + \mathbf{B}\boldsymbol\varepsilon\sqrt{\Delta t}
    \label{eq:coupled_OU}
\end{equation}
and setting $c_i(X)=1$.
Here $\mathbf{r} \in\mathbb{R}^K$ denotes the value of the stochastic process, $\boldsymbol\mu \in\mathbb{R}^K$ is the mean, $\mathbf{A}\in\mathbb{R}^{K\times K}$ defines the coupling between the different components, $\mathbf{B}\in\mathbb{R}^{K\times q}$ determines the noise scaling, and $\boldsymbol\varepsilon\in\mathbb{R}^q$ are $q$ independent Brownian motions,  that is, $\varepsilon_i(t)\sim N(0, 1)$. In what follows, we choose $q = K$.
The coupled OU process offers greater flexibility compared to the previously used independent scalar OU processes. Namely, interaction between the modes is encoded in $\mathbf{A}$, and incorporating the mean $\boldsymbol \mu$ removes the need to include the mean profile $\xi_0$ in the model.

The multiplicative noise model is presently defined by setting $c_i(X) = \left(a_i + b_i(\xi_i\cdot X)^2 \right)$ to modulate the noise with the energy of the solution expressed in terms of the POD modes. The parameters $a_i$ and $b_i$ will be tuned during the training procedure, and allow the learned model to incorporate state-dependent forcing in a straightforward manner.

Modeling the error of a nonlinear dynamical system with a linear stochastic process, such as the coupled OU process \eqref{eq:coupled_OU}, is a deliberate simplification. Although the stochastic process itself is linear, its influence on the forecast is nonlinear owing to its integration within the nonlinear system dynamics. This provides a dynamically consistent approach to incorporate the impact of unresolved scales \cite{wu2024learning}.
In particular, the additive coupled OU process provides a tractable and interpretable method for introducing temporally correlated noise. The multiplicative model represents a step towards a more realistic state-dependent parametrization.

In the present study, all model variables are obtained through an optimization procedure detailed in the next subsection. Specifically, no restrictions are imposed on the learnable parameters; the mean $\boldsymbol \mu$, the coupling matrix $\mathbf{A}$, the noise scaling $\mathbf{B}$, and the parameters $a_i$ and $b_i$ are determined without manual intervention. This approach is deemed feasible given the low dimensionality of the noise model, which results in only 160 tunable model parameters in the multiplicative case.

\subsubsection{Model training}\label{subsubsec:model_training}
The model parameters are obtained via an iterative optimization procedure. The high-fidelity data, spanning 500 MTU and stored at intervals of 0.005 MTU, are divided into calibration data and validation sets in an 80/20 split. Denoting the model parameters collectively by $\theta$, and the obtained ensemble prediction by $X(\theta)$, the loss is defined as 
\begin{equation}
    \text{Loss per batch} = \frac{1}{N_b}\frac{1}{N_t}\sum_{i=1}^{N_b}\sum_{j=1}^{N_t} \mathrm{CRPS} (X(\theta), C(q), t_j),
    \label{eq:loss_per_batch}
\end{equation}
where $C(q)$ is the coarsened reference configuration.
The batch size is denoted by $N_b$ and depends on the chosen trajectory length $N_t$, such that $N_b N_t$ maintains a comparable magnitude when $N_t$ is varied.
Furthermore, a batch comprises $N_b$ ensembles, each consisting of 20 ensemble members during training. The CRPS in \eqref{eq:loss_per_batch} is computed per $X$-variable separately and subsequently summed over all variables to obtain the loss per batch.

The model is trained during 400 epochs, after which the loss per batch was found to have stabilised. Each epoch consists of 125 batches. The Adam optimizer \cite{kingma2014adam} is used with a fixed learning rate set at $10^{-2}$ for training the additive noise model and $10^{-4}$ for the multiplicative noise model. A series of trajectory lengths is used, ranging from $N_t=2$ (0.01 MTU) to $N_t=200$ (1 MTU). A full description of the batch size, trajectory length, and final measured loss per batch is provided in Appendix \ref{app:details}.

\begin{figure}[h]
    \centering
    \includegraphics[width=0.49\linewidth]{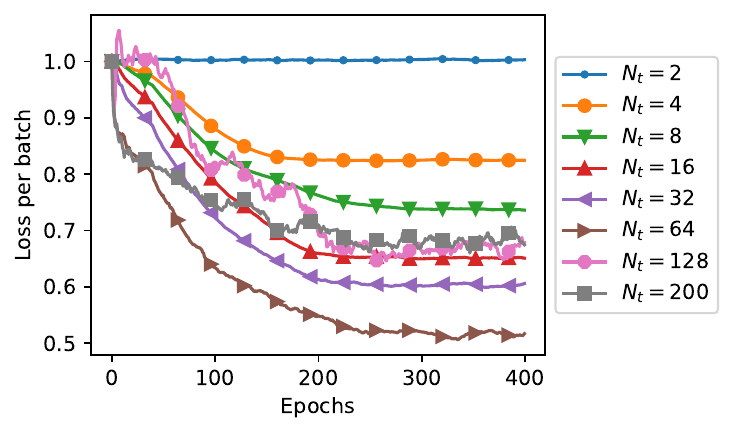}
    \includegraphics[width=0.49\linewidth]{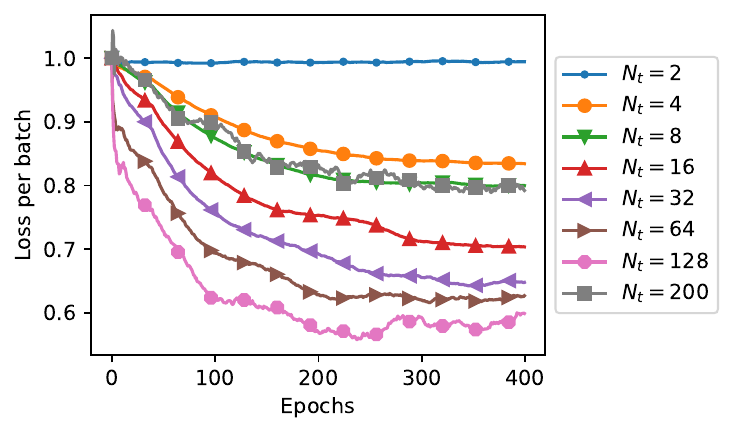}
    \includegraphics[width=0.49\linewidth]{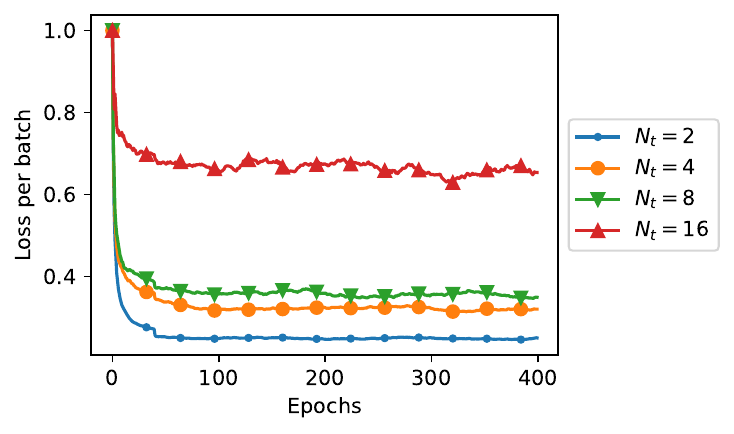}
    \includegraphics[width=0.49\linewidth]{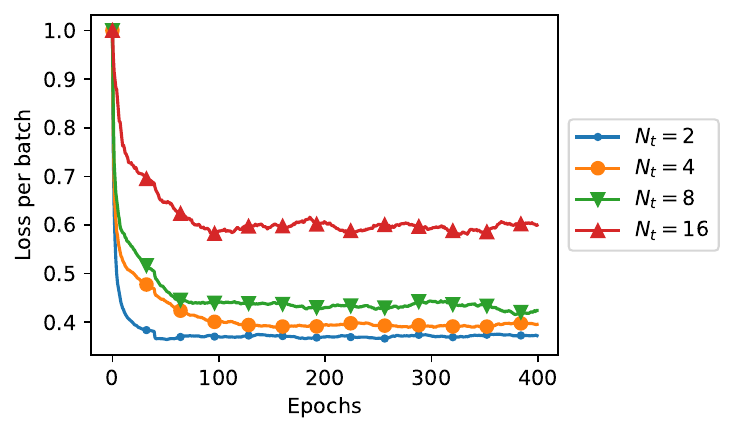}
    \caption{Evolution of the loss per batch over the training duration, for $c=4$ (left column) and $c=10$ (right column), for the additive noise model (top row) and the multiplicative noise model (bottom row). The displayed values are normalised by the initial value for visual comparison, and a moving average is applied for clarity.}
    \label{fig:loss_per_batch_epoch}
\end{figure}

The loss per batch during training is presented in Fig. \ref{fig:loss_per_batch_epoch} for different values of $N_t$. We observe that a short trajectory length for the additive noise model yields only a marginal decrease in the loss. Since each training trajectory starts from a single perfect initial condition, the error and spread after two time steps are small. Hence, for such short time intervals, the initial loss shows is barely improved upon during training. This effect diminishes for longer trajectories, and a decrease of the loss during training is realised until the trajectories become too long. Unstable training behavior is observed for $N_t=128$ and $N_t=200$, which is attributed to the chaotic dynamics of the system. Specifically, long model trajectories will at some point diverge from the reference trajectories. Minor changes in model parameters will induce significant changes in trajectory due to high system sensitivity, and too much effort is put in correcting long lead-time dynamics. This can lead to unstable gradients and inefficient optimization. This issue could be mitigated by weighting the loss function such that long lead-time predictions contribute less to the loss \cite{melchers2023comparison}; however, this is not further investigated here.

The evolution of the loss for the multiplicative noise model differs from that of the additive noise model. We observe that a short trajectory length ($N_t=2$) leads to a significant, rapid decrease of the loss, owing to the increased model parameters when compared to the additive noise model. Furthermore, training becomes less effective even at modest trajectory lengths, suggesting higher sensitivity to the model parameters as the trajectory length increases.

\section{Model performance}\label{sec:model_performance}
The prediction quality of each model is assessed at short lead times (i.e., `weather forecasts' of the L96 system) and at long lead times (i.e., the `climate' of the system).

The weather forecasts are carried out via ensemble predictions for 2 MTU, where each ensemble comprises 50 members. The slow variables sampled from 100 high-fidelity solution snapshots are used as initial conditions. These snapshots, each separated by 10 MTU, lie outside of the training/validation dataset.
We consider perfectly known initial conditions in Section \ref{subsec:short_time_predictions} and perturbed initial conditions in Section \ref{subsec:perturbed_short_time_predictions}. The prediction ensembles are assessed via the mean squared error (MSE), the spread-error relationship, and the CRPS. The predicted climate is assessed via a single long-time ensemble prediction, computed over 3000 MTU. The resulting distributions of the slow variables are compared via the Kolmogorov--Smirnov distance and the Hellinger distance. All error assessment metrics are further detailed in Appendix \ref{app:metrics}.
A summary of the models tested in this section is given in Table \ref{tab:model_configurations}.

\begin{table}[h]
    \centering
    \caption{Summary of the tested model configurations.}
    \begin{tabular}{l|l|l|l}
        Name & Model type & Model form & Noise  \\
        \hline
        poly\_gauss & Derivative fitting & Local & White, independent, additive\\
        poly\_ou    & Derivative fitting & Local & Correlated, independent, additive\\
        svd\_gauss  & Derivative fitting & Global & White, independent, additive\\
        svd\_ou     & Derivative fitting & Global & Correlated, independent, additive\\
        crps\_ou2   & Trajectory fitting, $N_t=2$ & Global & Correlated, coupled, additive\\
        crps\_ou8   & Trajectory fitting, $N_t=8$ & Global & Correlated, coupled, additive\\
        crps\_ou32  & Trajectory fitting, $N_t=32$ & Global & Correlated, coupled, additive\\
        crps\_mult2 & Trajectory fitting, $N_t=2$ & Global & Correlated, coupled, multiplicative\\
        crps\_mult8 & Trajectory fitting, $N_t=8$ & Global & Correlated, coupled, multiplicative\\
    \end{tabular}
    \label{tab:model_configurations}
\end{table}

\subsection{Short-time predictions} \label{subsec:short_time_predictions}
\begin{figure}[h]
    \centering
    \includegraphics[width=0.49\linewidth]{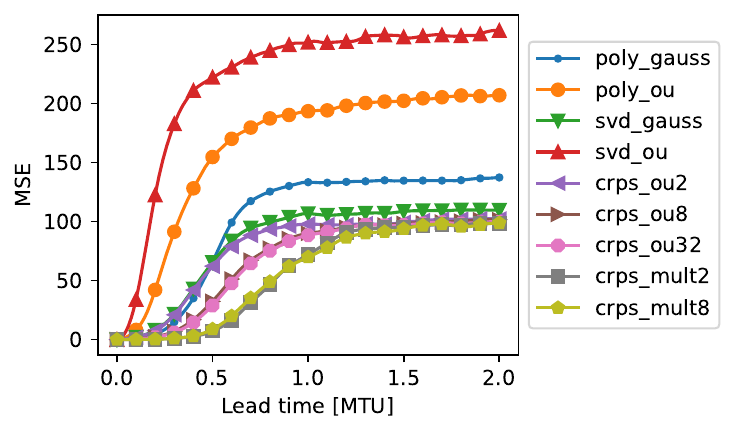}
    \includegraphics[width=0.49\linewidth]{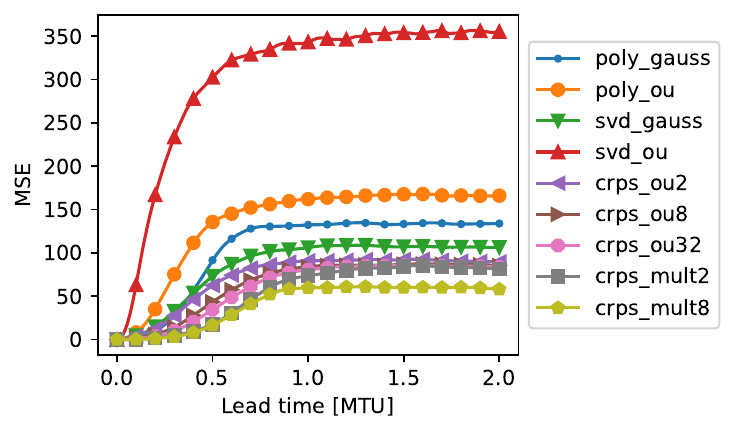}
    \caption{Comparison of the average MSE for short lead times, for $c=4$ (left) and $c=10$ (right). The MSE is measured per $X$-variable in ensembles consisting of 50 members, simulated from 100 different initial conditions.}
    \label{fig:MSE}
\end{figure}

The mean squared error (MSE) for short lead times is shown in Fig. \ref{fig:MSE}, for both test cases. Qualitatively, the results are consistent across both test cases. The models with independent correlated noise (`poly\_ou', `svd\_ou') are found to perform poorly. This is attributed to poor estimation of the correlation parameters in the autoregressive model. Namely, a visual comparison of the autocorrelation of the fitted autoregressive processes and the measured residuals and POD time series (Fig. \ref{fig:autocorrelation_models} in Appendix \ref{app:details}) suggests that the fitted temporal correlation parameters are too high. The uncorrelated derivative-fitting-based models (`poly\_gauss', `svd\_gauss') initially exhibit a similar growth of the MSE, but the global model saturates at a lower value than the local model, suggesting that incorporating spatial correlation patterns improves prediction accuracy. 

A further increase in average accuracy is observed for the CRPS-based models. Increasing the trajectory length generally leads to parametrizations with improved accuracy, up to a certain limit. For example, the additive noise model trained at $N_t=8$ is more accurate than that trained at $N_t=2$, and is approximately as accurate as that trained at $N_t=32$. The multiplicative models exhibit an even smaller MSE regardless of the trajectory length used in training, suggesting that the additional degrees of freedom allow for increased prediction accuracy. 

\begin{figure}[h]
    \centering
    \includegraphics[width=0.49\linewidth]{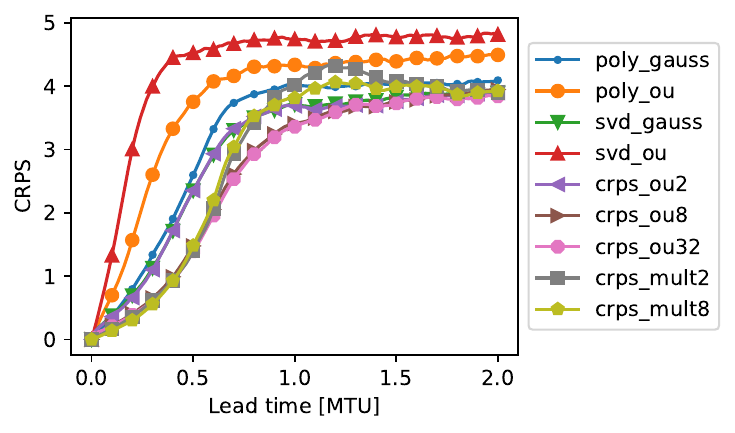}
    \includegraphics[width=0.49\linewidth]{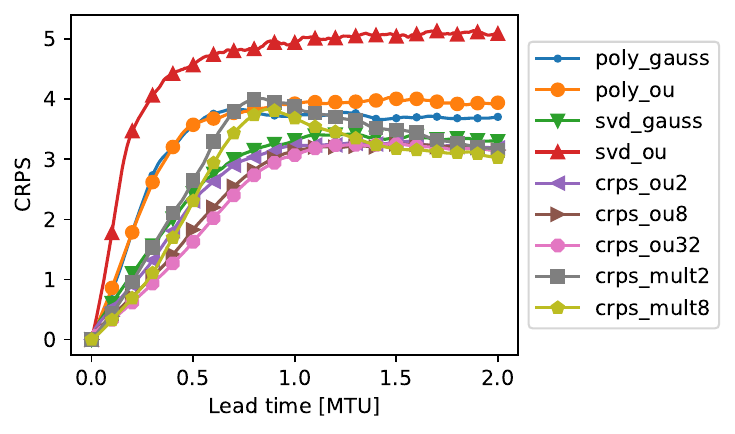}
    \caption{Comparison of the average CRPS for short lead times, for $c=4$ (left) and $c=10$ (right). The CRPS is measured per $X$-variable in ensembles consisting of 50 members, simulated from 100 different initial conditions.}
    \label{fig:CRPS}
\end{figure}

The average CRPS is shown in Fig. \ref{fig:CRPS}. The derivative-fitting-based models with correlated noise show rapid growth of the CRPS, owing to the large mean error of these ensembles. However, this becomes less pronounced when $c=10$, where we observe that the local models (`poly\_gauss', `poly\_ou') perform approximately equally well. The global model with uncorrelated noise exhibits similar evolution of the CRPS as the additive noise model trained at $N_t=2$, and generally performs well. Overall, the lowest CRPS is obtained for the additive noise models trained with $N_t=8$ and $N_t=32$. For $c=4$, the multiplicative noise models show a similar evolution to these additive noise models. In contrast, for $c=10$, the multiplicative noise models yield a higher CRPS on average.

\begin{figure}[h]
    \centering
    \includegraphics[width=0.85\linewidth]{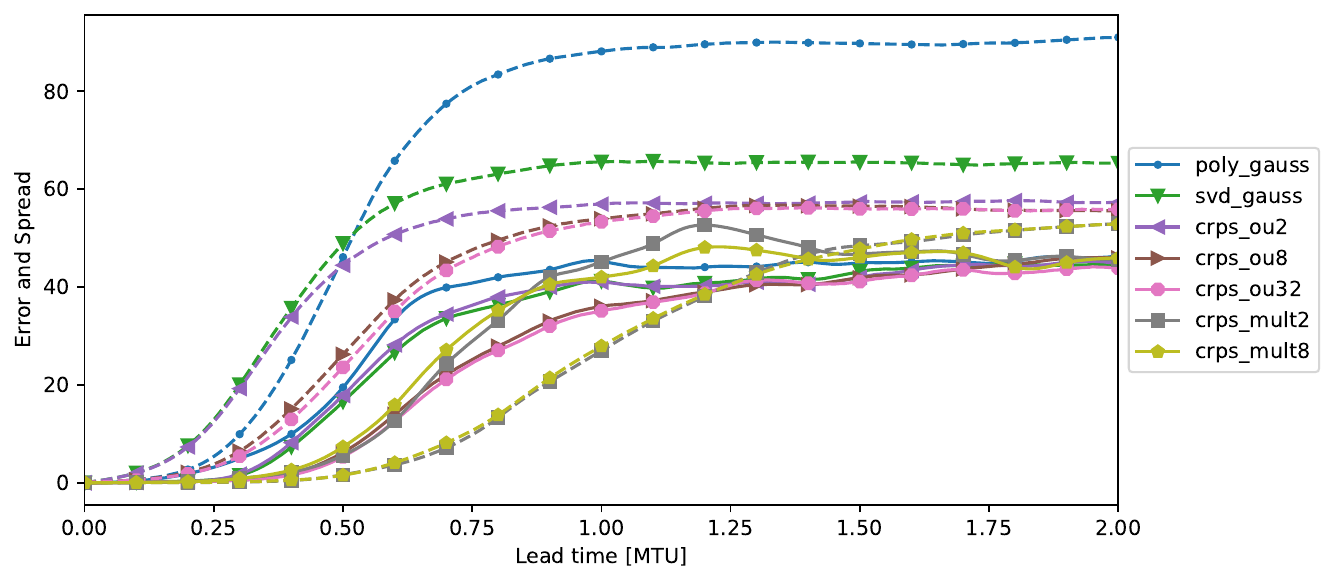}
    \includegraphics[width=0.85\linewidth]{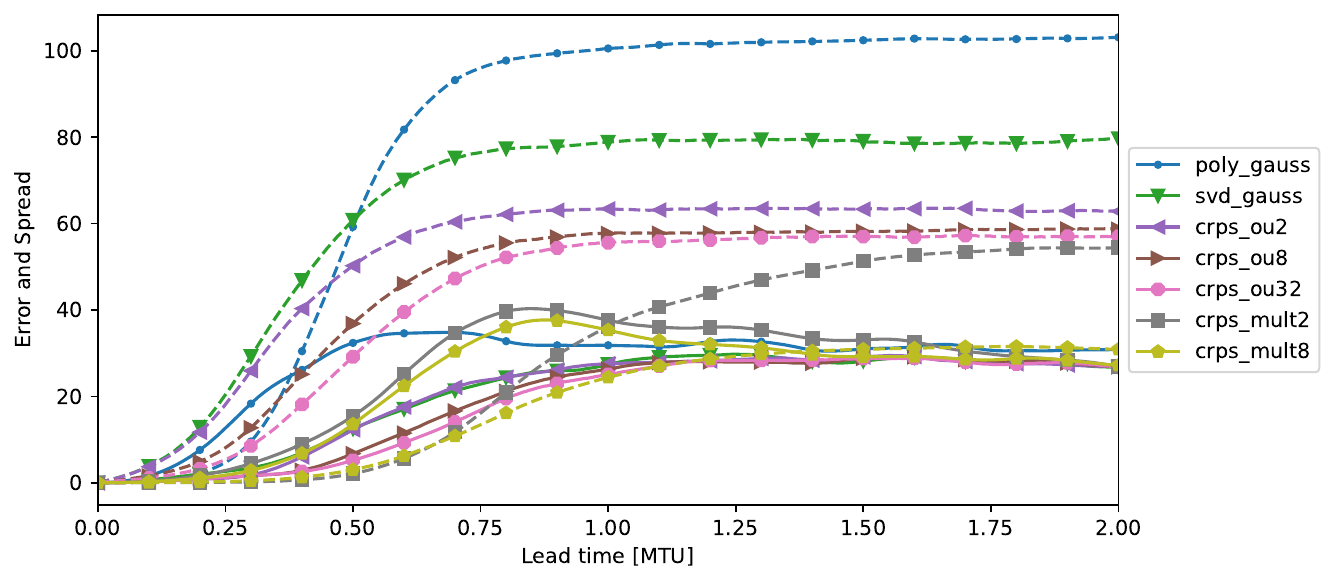}
    \caption{Comparison of ensemble error (solid lines) and ensemble spread (dashed lines) for short lead times using perfect initial conditions. The results of the `poly\_ou' and `svd\_ou' model configurations are not shown here due to their high errors. The error and spread are computed per $X$-variable in ensembles consisting of 50 members, simulated from 100 different initial conditions.}
    \label{fig:MSE_spread}
\end{figure}

Although the multiplicative models display a smaller MSE (Fig. \ref{fig:MSE}) and smaller spread (Fig. \ref{fig:MSE_spread}) than the other methods, this is not clearly reflected in the CRPS. The former metrics are measured as 2-norms, whereas the CRPS can be regarded as a measure of error and spread in the 1-norm. This suggests that the error and spread of the multiplicative model contains numerous small values that contribute little to the 2-norm but accumulate in their contribution to the 1-norm. In turn, this implies that the corresponding predictions are both accurate and sharp, in line with the training objective specified by the CRPS loss function.

A visual comparison between the ensemble error and spread is provided in Fig. \ref{fig:MSE_spread}. A necessary condition for statistical consistency in ensembles is that the mean ensemble variance estimate equals the squared ensemble mean error \cite{leutbecher2009diagnosis}. Here, the ensemble error is computed following \eqref{eq:MSE_spread}. We observe that all additive noise models are overdispersive, since the ensemble spread is larger than the ensemble error, which is maintained for all lead times. The forecast ensembles obtained with the multiplicative model are found to be underdispersed, implying an overconfident forecast based on the measured error. This is maintained until a lead time of approximately 1.5 MTU for $c=4$ and 1 MTU for $c=10$. 

For the trained models, it could be desirable to diminish the weight of the self-spread term in the CRPS loss to possibly favor statistical consistency. This adaptation requires an additional hyperparameter or a term penalizing the discrepancy between the ensemble spread and MSE. Such a modification might diminish the ensemble sharpness during training without sacrificing accuracy. We investigate the use of a hyperparameter to scale the self-spread term in Section \ref{subsec:self_spread}.

\subsection{Long-time predictions}\label{subsec:long_time_predictions}
\begin{figure}[h]
    \centering
    \includegraphics[width=0.49\linewidth]{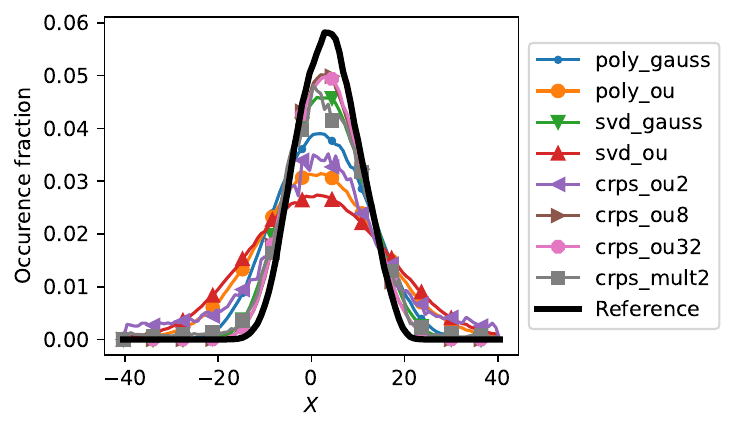}
    \includegraphics[width=0.49\linewidth]{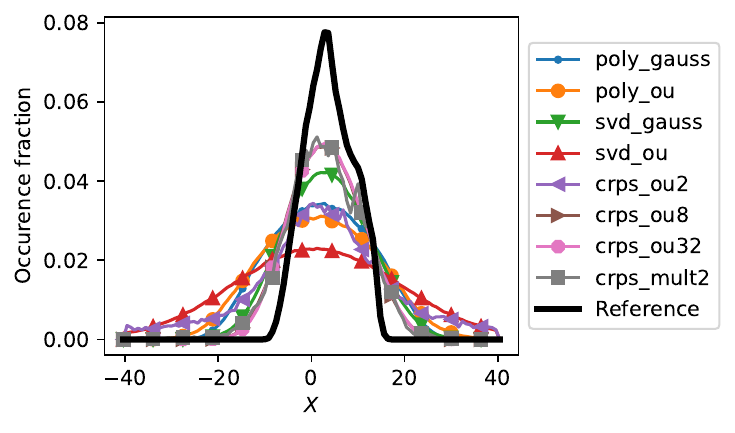}
    \caption{Comparison of the climate of the L96 system for $c=4$ (left panel) and $c=10$ (right panel), defined as the PDFs of the $X$-variables.}
    \label{fig:long_time_pdfs_corr}
\end{figure}

The climate prediction using the stochastic parametrizations is shown in Fig. \ref{fig:long_time_pdfs_corr}, defined via the empirical PDFs of the slow variables over a long simulation time (3000 MTU). Here, the `crps\_mult8' is not included in the results, as it exhibited numerical instabilities during extended simulation. Generally, all stochastic parametrizations produce qualitatively similar symmetric PDFs with an excessively large standard deviation. An improvement of the CRPS-based models is found compared to the derivative-fitting-based models, which becomes clearer when increasing the trajectory length. This is summarised in Table \ref{tab:long_time_distances}, indicating that the additive noise models learned using longer trajectory lengths outperform the other tested models. However, the results are not as accurate as those presented in literature (see, e.g., \cite{gagne2020machine}, \cite{crommelin2008subgrid}), suggesting that combining CRPS-based learning with other state-dependent model forms might be beneficial for long-time prediction accuracy.

\begin{table}[h]
    \centering
    \caption{Kolmogorov--Smirnov (KS) and Hellinger ($\mathrm{H}^2$) distances, computed from the empirical PDFs of the $X$-variables over long simulation times. The PDFs are shown in Fig. \ref{fig:long_time_pdfs_corr}.}
    \begin{tabular}{l|c|c|c|c}
    & \multicolumn{2}{c|}{$c=4$} & \multicolumn{2}{c}{$c=10$} \\
    \cline{2-5}
    Model & $\mathrm{KS}$ & $\mathrm{H}^2$ & $\mathrm{KS}$ & $\mathrm{H}^2$ \\
    \hline
      poly\_gauss  & 0.136          & 0.258         & 0.238         & 0.466 \\
      poly\_ou     & 0.204          & 0.390         & 0.269         & 0.513  \\
      svd\_gauss   & 0.0693         & 0.160         & 0.160         & 0.377 \\
      svd\_ou      & 0.240          & 0.449         & 0.330         & 0.624 \\
      crps\_ou2    & 0.173          & 0.383         & 0.241         & 0.515  \\
      crps\_ou8    & \textbf{0.0504}& 0.109         & 0.129         & 0.303  \\
      crps\_ou32   & 0.0506         & \textbf{0.108}& \textbf{0.122}& \textbf{0.299}  \\
      crps\_mult2  & 0.0838         & 0.199         & 0.138         & 0.313
    \end{tabular}
    \label{tab:long_time_distances}
\end{table}

    %

\subsection{Short-time predictions with perturbed initial conditions}\label{subsec:perturbed_short_time_predictions}
The ensemble predictions at short lead times presented in Section \ref{subsec:short_time_predictions} are repeated with perturbed initial conditions to assess the robustness of the stochastic parametrizations. The slow variables constituting the initial conditions for each ensemble are perturbed by a random amount, achieved by sampling $r\sim\mathcal{N}(0, 0.1\hspace{0.3mm}S^2)$. Here $S^2$ is the sample variance of the slow variables in the reference data, measured as $S^2=6.45$ for $c=4$ and $S^2=5.08$ for $c=10$. The perturbations are applied independently per variable and per initial condition, however, they are the same for each model configuration to ensure a fair comparison. 

\begin{figure}[h]
    \centering
    \includegraphics[width=0.85\linewidth]{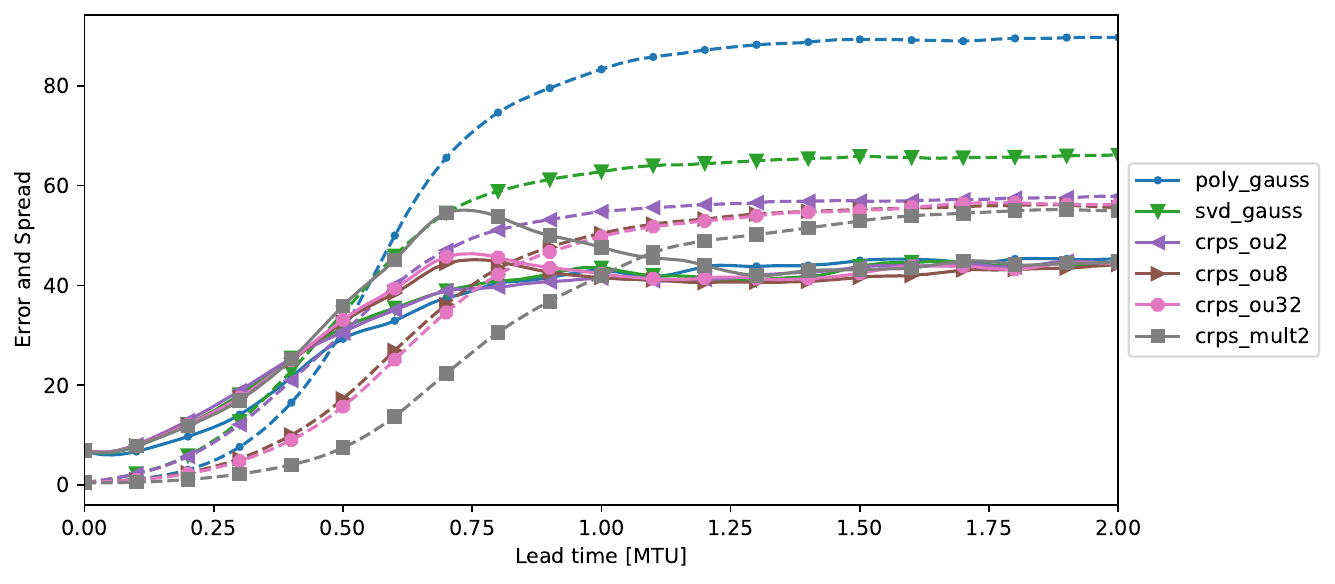}
    \includegraphics[width=0.85\linewidth]{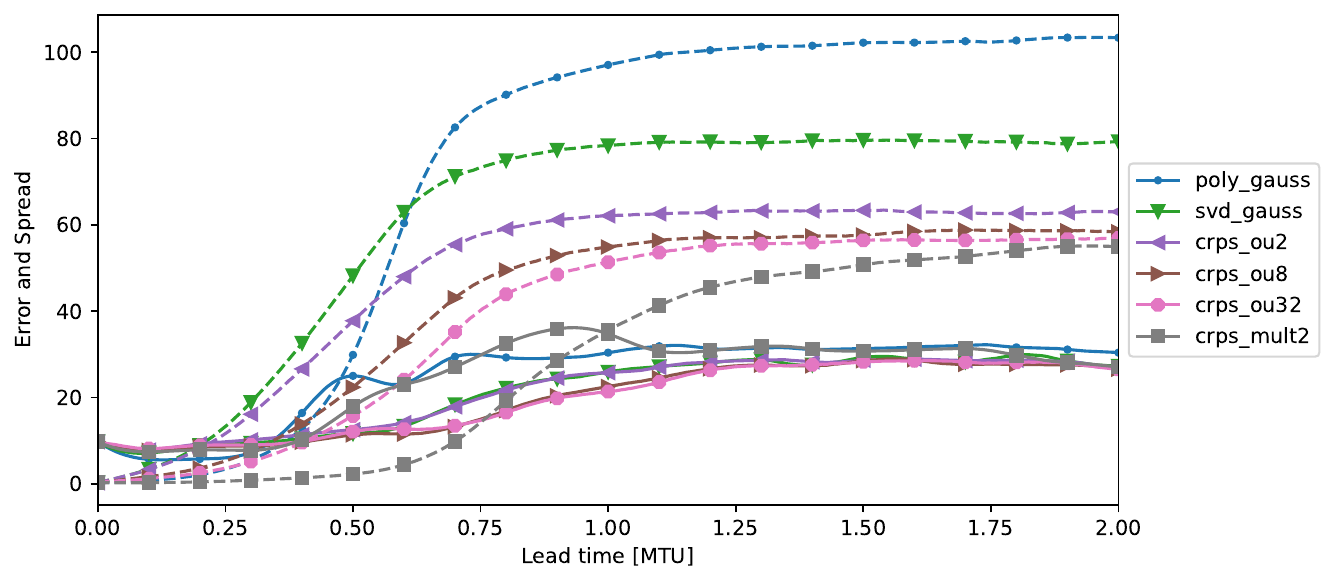}
    \caption{Comparison of ensemble error (solid lines) and ensemble spread (dashed lines) for short lead times using perturbed initial conditions. The results of the `poly\_ou' and `svd\_ou' model configurations are not shown here due to their high errors. The error and spread are computed per $X$-variable in ensembles consisting of 50 members, simulated from 100 different initial conditions.}
    \label{fig:MSE_spread_pert}
\end{figure}

The resulting ensemble error and spread are shown in Fig. \ref{fig:MSE_spread_pert}. The resulting spread values are qualitatively similar to those from the perfect initial condition experiments (Fig. \ref{fig:MSE_spread}). The multiplicative noise model exhibits an increased error growth until approximately 0.75 MTU when compared to the other models. This increased error is attributed to the state-dependent forcing losing accuracy when perturbations are introduced to the initial conditions. Overall, the qualitative similarity observed in the error of the additive models suggests that the learned parametrizations are robust against initial condition uncertainty.

\subsection{Sensitivity to self-spread penalty}\label{subsec:self_spread}
The previous subsections highlight that the multiplicative noise models generate sharp ensemble forecasts, particularly visible at shorter lead times. This is a consequence of using the CRPS as the loss function, which balances accuracy and ensemble spread, favouring slightly dispersed forecasts that still reproduce the observed value during training.
From a statistical perspective, however, forecast reliability requires that the ensemble spread equals the ensemble error \cite{talagrand1999evaluation, leutbecher2009diagnosis}. We therefore investigate the role of the spread term in the CRPS loss function and quantify its influence on short-term forecasts.

To this end, we modify the loss function \eqref{eq:CRPS_loss} so that the self-spread penalty is scaled by a factor $\alpha$, \begin{equation} \label{eq:CRPS_scaled}
    \begin{split}
        \mathrm{Loss}(\theta, q, \alpha) &= \sum_{j=1}^{N_t}\mathrm{CRPS}(z_\theta, C(q), t_j, \alpha), \\
        \mathrm{CRPS}(z_\theta, C(q), t_j, \alpha) &= \frac{1}{M}\sum_{i=1}^M \left|z_\theta^i(t_j) - C(q(t_j)) \right| - \frac{\alpha}{2M^2} \sum_{j=1}^M\sum_{k=1}^M\left| z_\theta^i(t_j) - z_\theta^k(t_j) \right|, 
    \end{split}
\end{equation}
where $\alpha \in[0, 1]$. We henceforth refer to \eqref{eq:CRPS_scaled} as the \textit{scaled CRPS loss function}. The addition of $\alpha$ allows for training models that minimize the ensemble mean error while reducing the importance of the ensemble spread. 
The standard CRPS loss function \eqref{eq:CRPS_loss} is recovered when $\alpha=1$. The forecast spread term is progressively diminished when $\alpha$ decreases, and vanishes entirely when $\alpha=0$, leaving only the ensemble mean error in the loss function. The value of $\alpha$ cannot exceed 1 to ensure nonnegativity of the scaled CRPS.

The effects of scaling the CRPS loss function are assessed by repeating the model training and forecasting procedures for a sequence of values of $\alpha$. A single model form is selected for the sake of clarity. Here, we adopt the multiplicative noise model and use a trajectory length $N_t=8$ to train the model, corresponding to the `crps\_mult8' setting considered in the previous subsections. This model was found to produce accurate yet sharp ensemble forecasts, and is therefore suited to demonstrate the ability to change the ensemble spread via choice of the loss function.

\begin{figure}[h]
    \centering
    \includegraphics[width=0.49\linewidth]{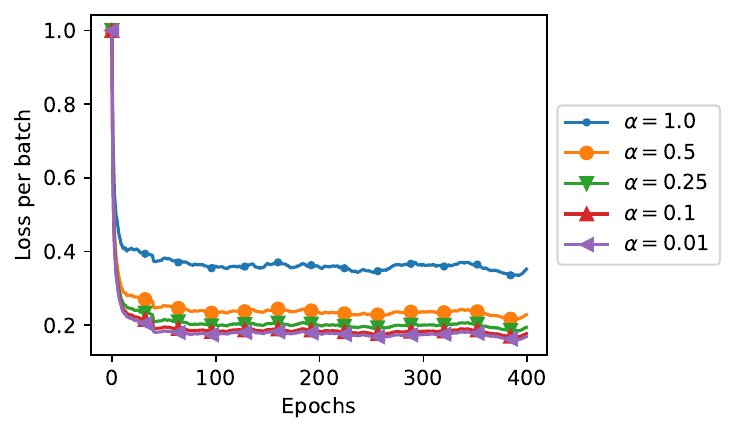}
    \includegraphics[width=0.49\linewidth]{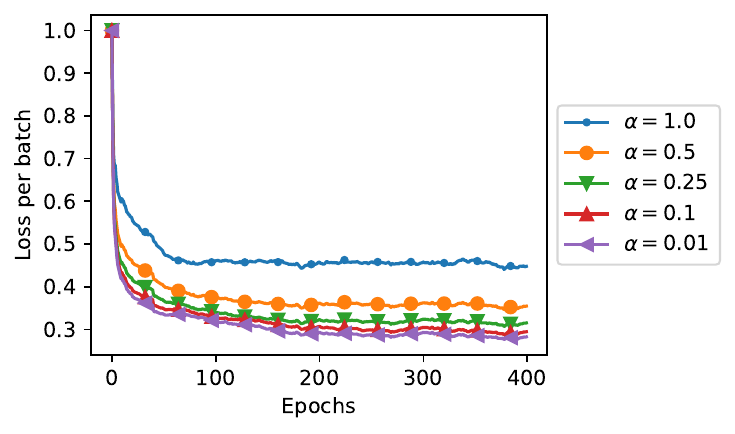}
    \caption{Evolution of the loss per batch over the training duration, for $c=4$ (left panel) and $c=10$ (right panel), for several values of $\alpha$ in the scaled CRPS loss function \eqref{eq:CRPS_scaled}. The displayed values are normalised by the initial value for visual comparison, and a moving average is applied for clarity.}
    \label{fig:loss_per_batch_epoch_scaled}
\end{figure}

The evolution of the training loss is shown in Fig. \ref{fig:loss_per_batch_epoch_scaled}, demonstrating consistent behaviour across the two test cases. Starting from the same initial loss value, the loss decreases at a similar rate for each values of $\alpha$. Smaller final losses are achieved for smaller $\alpha$. This is expected, since a reduced weight on the self-spread term in the scaled CRPS imposes fewer or less stringent constraints during training, thereby facilitating easier optimisation of the loss function.

\begin{figure}[h]
    \centering
    \includegraphics[width=0.85\linewidth]{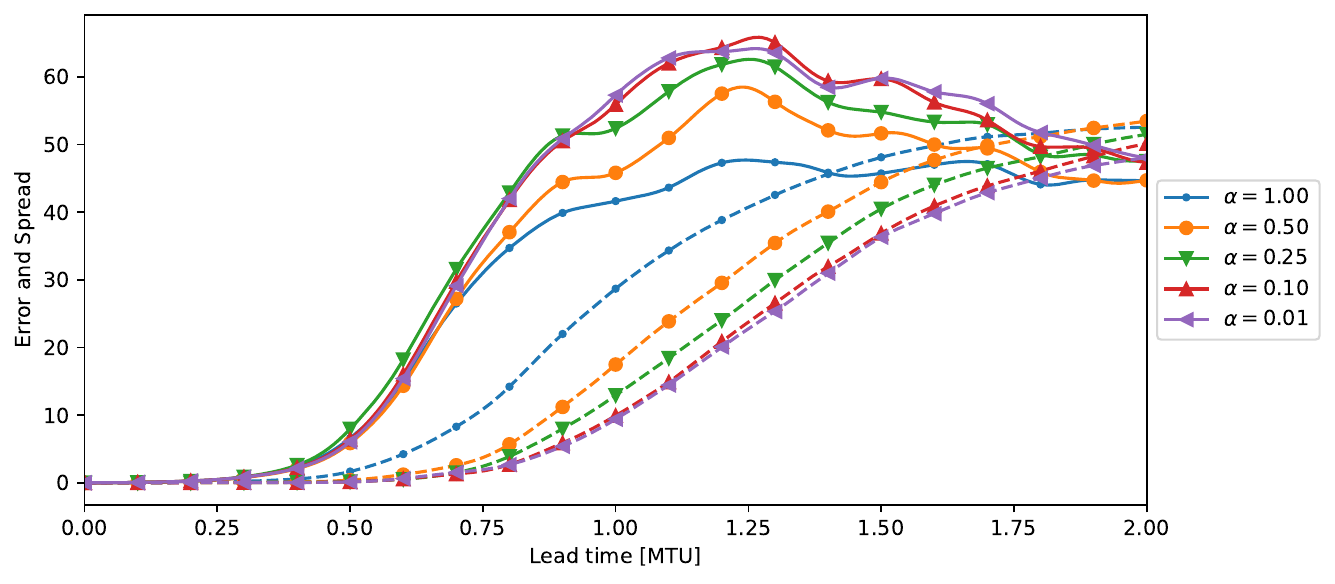}
    \includegraphics[width=0.85\linewidth]{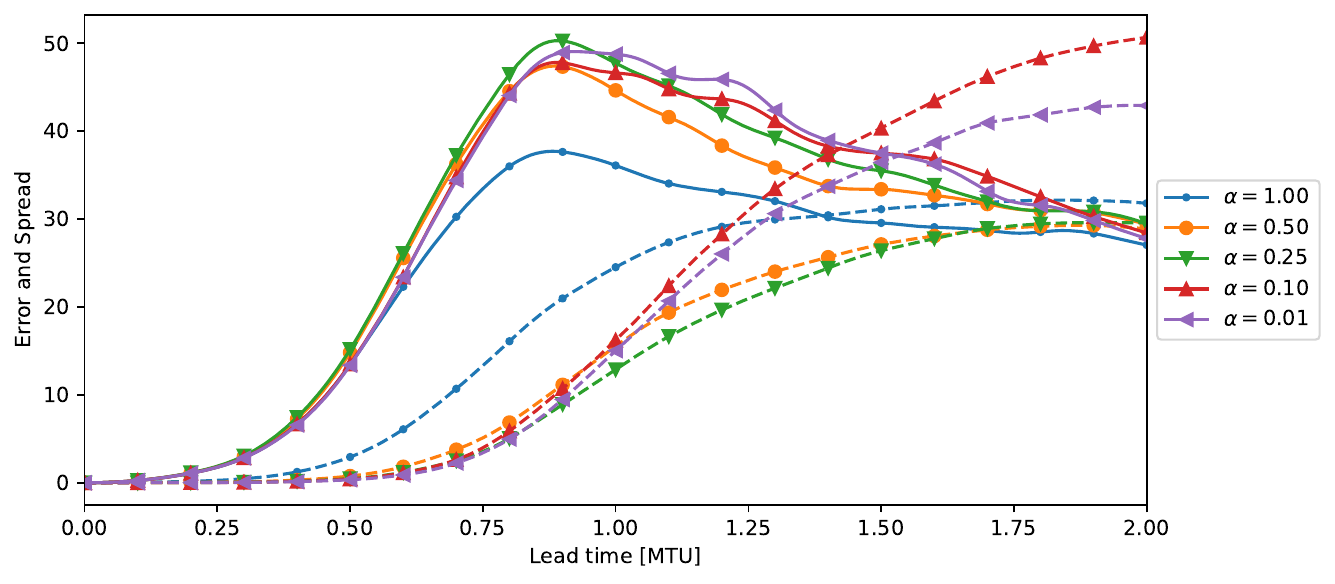}
    \caption{Comparison of ensemble error (solid lines) and ensemble spread (dashed lines) for short lead times, for several values of $\alpha$ in the scaled CRPS loss function \eqref{eq:CRPS_scaled}. The error and spread are computed per $X$-variable in ensembles consisting of 50 members, simulated from 100 different initial conditions.}
    \label{fig:MSE_spread_scaled}
\end{figure}

A comparison between the ensemble error and ensemble spread is provided in Fig. \ref{fig:MSE_spread_scaled}. A clear effect on the ensemble spread is observed in both test cases, with the forecast sharpness increasing as $\alpha$ decreases. This occurs because, when $\alpha=1$, the CRPS formulation favours a larger spread to minimise the loss. At the same time, the ensemble error increases with decreasing $\alpha$. This suggests that retaining the self-spread term in the CRPS is beneficial for reducing mean error, despite not directly contributing to the error itself. Larger values of $\alpha$ regularise training by preventing a collapse of ensemble members, thereby reducing overfitting and generally yielding more stable predictions.

For the same reason, the spread observed for $c=10$ and small values of $\alpha$ initially remains small but later exceeds the spread obtained with $\alpha=1$. This behaviour is also attributed to the regularising property of the self-spread term in the loss function. Namely, models trained with small $\alpha$ are more sensitive to perturbations, and their forecast spread is therefore amplified more at later simulation times. These results demonstrate that the ensemble spread and model robustness can be systematically controlled through the choice of loss function, specifically via the weighting of the self-spread term.

\section{Conclusion}\label{sec:conclusion}
In this paper, we investigated the calibration of stochastic parametrizations for ensemble forecasting using trajectory learning based on the continuous ranked probability score (CRPS). Using the two-scale Lorenz '96 system as a testbed, we demonstrated that this approach yields models that provide both accurate and sharp short-term forecasts. Both additive and multiplicative noise models, leveraging global spatial patterns from POD modes, were found to outperform models based on derivative fitting. The improved accuracy at short lead times suggests that CRPS-based trajectory fitting is promising approach when combined with data assimilation, where accurate short-term forecasts are essential. Moreover, the self-spread term in the CRPS was found to have a stabilising effect on the model. Ultimately, this work establishes CRPS-based trajectory learning as a promising methodology for developing data-driven stochastic parametrizations for complex dynamical systems, representing a first step towards its application in geophysical flows.

The stochastic models presented here serve as a proof-of-concept that trajectory learning can be applied to yield accurate ensemble forecasts. A natural next step is to develop additional state-dependent parametrizations using tools such as neural networks, which could capture more complex nonlinear error dynamics. However, many challenges remain. As previously observed in deterministic trajectory fitting, the length of the training trajectory remains a critical hyperparameter. Performance generally improved with longer trajectories, but became unstable at certain lengths due to the chaotic system dynamics. For complex parametrizations, such as those based on neural networks, this sensitivity might become even more pronounced. Future work could explore strategies to mitigate these adverse effects, for instance, by assigning greater weight on short lead times within the loss function to emphasize initial accuracy.

\subsubsection*{Acknowledgements}
SE gratefully acknowledges the Royal Society of Arts and Sciences in Gothenburg (KVVS) for funding the research visit to Imperial College London that started this work. This work was supported in part by the Swedish Research Council (VR) through Grant No. 2022-03453 and by the European Research Council (ERC) Synergy grant “Stochastic Transport in Upper Ocean Dynamics” (STUOD) – DLV-85640.

\bibliographystyle{abbrv}
\bibliography{refs}

\appendix

\section{Model and training details} \label{app:details}
This appendix provides supplementary details for the models presented in the text.

\begin{figure}[h]
    \centering
    \includegraphics[width=0.41\linewidth]{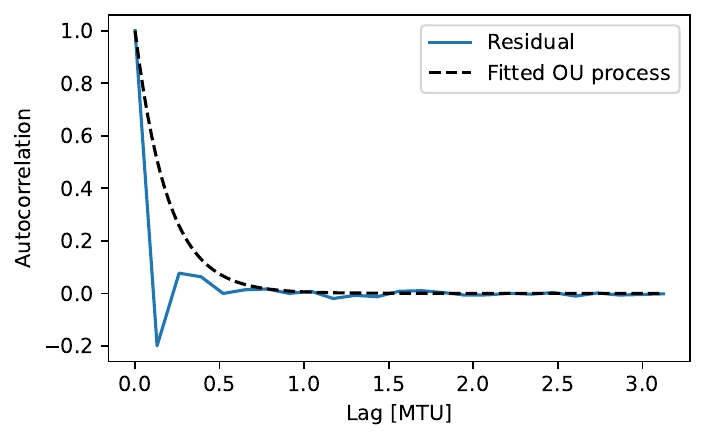}
    \includegraphics[width=0.41\linewidth]{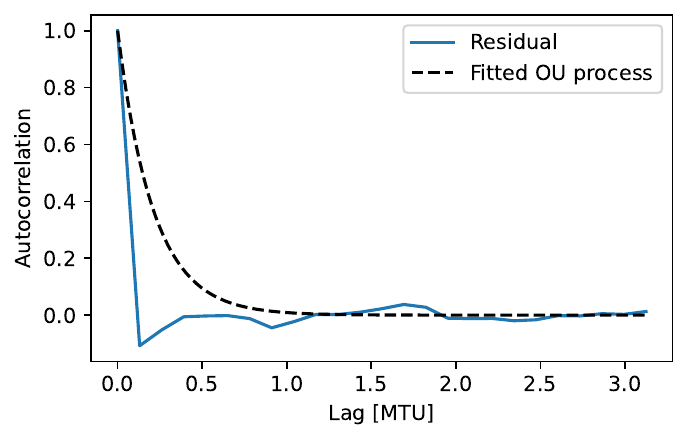}
    \includegraphics[width=0.41\linewidth]{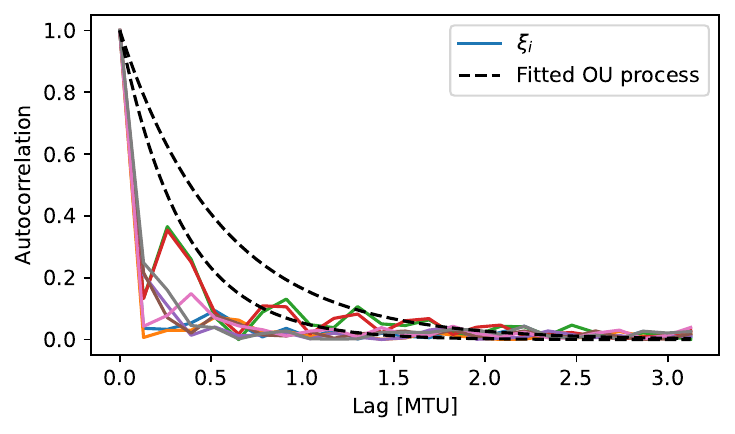}
    \includegraphics[width=0.41\linewidth]{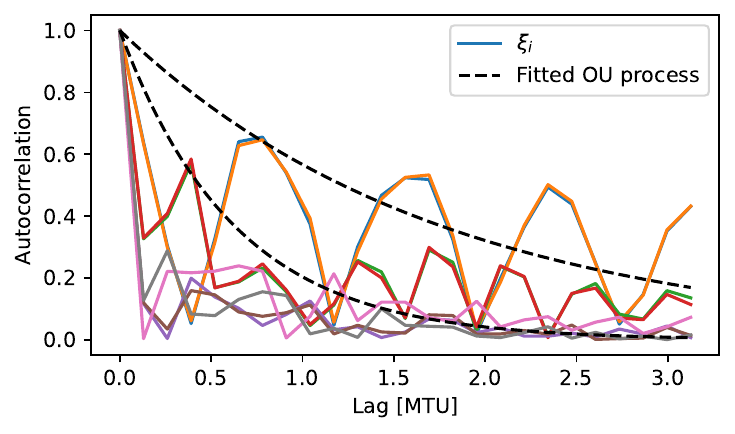}
    \caption{Comparison of sampled autocorrelation and autocorrelation of fitted OU processes, shown for the residuals of the polynomial fit (top row) and the POD modes (bottom row), for $c=4$ (left column) and $c=10$ (right column). The absolute value of each correlation is shown for the POD modes for comparison purposes.
    For the POD modes, the measured autocorrelation is shown using solid lines, and only the autocorrelation of the highest and lowest autoregressive coefficient are displayed.}
    \label{fig:autocorrelation_models}
\end{figure}

The parameters for the derivative-fitting-based models with correlated noise, introduced in Section \ref{subsec:derivative_fitting_models}, are summarised in Table \ref{tab:autocorrelation_parameters_pod}. These models use increments generated by an Ornstein--Uhlenbeck process, discretised as a first-order autoregressive (AR(1)) model, for the stochastic model component. The specific autocorrelation coefficients per POD mode are provided. Fig. \ref{fig:autocorrelation_models} provides a visual comparison between the measured autocorrelation from the data and the autocorrelation function of the AR(1)-processes with fitted parameters. As discussed in Section \ref{subsec:short_time_predictions}, this figure highlights a potential overestimation of the temporal correlation.

The training parameters used for the CRPS-based models, as described in Section \ref{subsubsec:model_training}, are given in Table \ref{tab:length_batch_loss}. A series of trajectory lengths $N_t$ were used to investigate the effect on model performance. The corresponding batch sizes $N_b$ are chosen to maintain a comparable computational load per epoch. For each value of $N_t$, the final loss per batch achieved during training for the additive and multiplicative models is included.

\begin{table}[h]
    \centering
    \caption{Parameters for the AR(1)-process used in the global model, per POD mode and per test case. }
    \begin{tabular}{c||c|c}
         & $c=4$ &  $c=10$\\
         \hline
        $\xi_1$ & 0.990 & 0.997 \\
        $\xi_2$ & 0.991 & 0.997 \\
        $\xi_3$ & 0.990 & 0.996 \\
        $\xi_4$ & 0.991 & 0.995 \\
        $\xi_5$ & 0.988 & 0.993 \\
        $\xi_6$ & 0.988 & 0.992 \\
        $\xi_7$ & 0.989 & 0.992 \\
        $\xi_8$ & 0.985 & 0.993 \\
    \end{tabular}
    \label{tab:autocorrelation_parameters_pod}
\end{table}
\begin{table}[h]
    \centering
    \caption{Overview of the training data for the CRPS-based models. Shown are the trajectory length $N_t$ and corresponding time interval in MTU, the batch size $N_b$, and the final losses per batch observed during training for the two considered test cases.}
    \begin{tabular}{l|c|c|c||c|c}
     & \multicolumn{3}{c||}{Training parameters} & \multicolumn{2}{c}{Final loss} \\
    \cline{2-6}
       & $N_t$ & MTU & $N_b$ & $c=4$ & $c=10$\\
        \cline{1-6}
     Additive noise  & 2 & 0.010 & 200 & 0.290 & 0.402\\
       & 4 & 0.020 & 100 & 0.492 & 0.694\\
       & 8 & 0.040 & 50  & 0.721 & 1.08\\
       & 16 & 0.080 & 25 & 1.01  & 1.52\\
       & 32 & 0.16 & 12 & 1.55  & 2.32\\
       & 64 & 0.32 & 6  & 2.97  & 4.55\\
       & 128 & 0.64 & 3 & 5.07  & 8.51\\
       & 200 & 1.0 & 2 & 14.3  & 14.3\\
        \cline{1-6}
     Multiplicative noise  & 2 & 0.010 & 200 & 0.0367  & 0.0707\\
       & 4 & 0.020 & 100 & 0.101  & 0.195\\
       & 8 & 0.040 & 50 & 0.275  & 0.410\\
       & 16 & 0.080 & 25 & 0.562  & 0.730\\
    \end{tabular}
    \label{tab:length_batch_loss}
\end{table}
\newpage

\section{Forecast assessment metrics}\label{app:metrics}
In the following definitions, we denote the truth by $y$ and denote by $x$ an ensemble forecast of size $M$, consisting of ensemble members $x_i,~i=1,\ldots,M$. All $x_j$ and $y$ are regarded as independent identically distributed random variables with mean $\mu$ and variance $\sigma^2$.

\subsection{Short lead time diagnostics}
\begin{enumerate}
    \item A general diagnostic of model accuracy is the mean square error (MSE), \begin{equation}
        \mathrm{MSE}(x,y, t) = \frac{\sum_{i=1}^M \left(x_i(t) - y(t)\right)^2}{M}.
        \label{eq:MSE}
    \end{equation}
    If many ensemble predictions are made, for example by using many different initial conditions, one may instead consider taking the mean value of the MSE.
    
    \item The spread-error relationship assesses the reliability of the forecast ensemble \cite{talagrand1999evaluation, leutbecher2009diagnosis}. A statistically consistent ensemble requires the ensemble spread ensemble and squared mean error to be equal. This is based on the identities of the ensemble mean error and the ensemble variance,
    \begin{align}
        \E\left(\frac{1}{M}\sum_{i=1}^M x_i-y\right)^2 &= \left(1+\frac{1}{M}\right)\sigma^2, \\
        \E\frac{1}{M}\sum_{i=1}^M\left(x_i-\frac{1}{M}\sum_{j=1}^M x_j\right)^2 &= \left(1-\frac{1}{M}\right)\sigma^2.
        \label{eq:MSE_spread}
    \end{align}
    This metric should ideally be applied after correcting the bias of the forecast, although this is not always feasible in practice. This relation can still be used when the mean forecast error is small compared to the random component of the forecast error.

    \item The Continuous Ranked Probability Score (CRPS), is a proper probabilistic forecast score used to assess the quality of a ensemble forecast. Given a empirical Cumulative Distribution Function (CDF) $F$, one can define the CRPS with respect to the observation $y_\mathrm{obs}$ by integrating the CDF as follows
\begin{align}
\operatorname{CRPS}(F, y_\mathrm{obs})=\int_{-\infty}^{\infty}(F(s)-\mathbb{I}_{\{s \geq y_\mathrm{obs}\}})^2 d s.
\end{align}
Here, lower scores indicate better forecast skill. For a discrete forecast, $\lbrace x_{1},...,x_{M}
\rbrace$ the CRPS score is approximated as follows, 
\begin{align}
\operatorname{CRPS}(x,y_\mathrm{obs})=\frac{1}{M} \sum_{i=1}^M\left|x_i-y_\mathrm{obs}\right|-\frac{1}{2 M^2} \sum_{i=1}^M \sum_{j=1}^M\left|x_i-x_j\right| .
\end{align}
Note that $1/2M^2$ can be replaced by $1/2M(M-1)$ when used as a estimator. The above CRPS definition is for a single point, and typically is subsequently averaged over points in which there is observations/data, either in time or space. The first term in the CRPS describes the 
accuracy to the data, and the second term penalizes self-spread as to enforce a sharper forecast. 
\end{enumerate}

\subsection{Long lead time diagnostics}
Diagnostics for long-time simulations generally concern the distribution of predicted variables and statistics of measured time series. Following We denote by $P$ the empirical forecast cumulative distribution function (CDF) and by $P'$ the verification (truth) CDF. The empirical probability distribution functions (PDFs) are respectively denoted by $p$ and $p'$. The PDFs are approximated by histograms with $M$ bins, where $p_i$ ($p'_i$) denotes the $i^\mathrm{th}$ bin of $p$ ($p'$). We compare the long-time distributions via the following diagnostics. 
\begin{enumerate}
    \item The Kolmogorov--Smirnov (KS) distance \begin{equation}
        \mathrm{KS}(P, P') = \max_i|P_i-P'_i|,
    \end{equation}
    which measures the largest mismatch between two distributions based on their CDFs. It is sensitive to shifts of distributions or spikes, since these yield large local discrepancies.

    \item The Hellinger distance
    \begin{equation}
        \mathrm{H}^2(p, p') = \frac{1}{\sqrt{2}}\sqrt{\sum_{i=1}^M \left(\sqrt{p_i} - \sqrt{p'_i}\right)^2},
    \end{equation}
    which is suited to measure the overlap between distributions defined by their probability density. Compared to the KS distance, it is more sensitive to shape differences in the distributions.
\end{enumerate}

\end{document}